\documentclass[a4paper, 12pt]{amsart}
\usepackage{mathrsfs}
\usepackage{amsmath}
\usepackage{amssymb}
\usepackage{verbatim}
\usepackage{xcolor}
\usepackage{CJK}
\usepackage{graphicx}
\usepackage[notref,notcite]{showkeys}
\DeclareGraphicsExtensions{.eps,.ps,.jpg,.bmp}

\theoremstyle{plain}

\newtheorem{theorem}{Theorem}[section]
\newtheorem{proposition}{Proposition}[section]
\newtheorem{corollary}{Corollary}[section]

\newtheorem{lemma}{Lemma}[section]

\theoremstyle{definition}
\newtheorem{definition}{Definition}[section]

\theoremstyle{remark}
\newtheorem{remark}{Remark}[section]
\newtheorem*{acknowledgment}{Acknowledgment}
\title{Symmetrization of plurisubharmonic functions on the Fano manifolds}
\author{Jingcao Wu$^1$}
\thanks{$^1$~Partially supported by CSC}

\begin{document}

\begin{abstract}
Given a compact complex manifold $Y$ with a negative line bundle $L\rightarrow Y$, we study the Schwarz-type symmetrization on the total space of $L$. We prove that this symmetrization does not increase the Monge-Amp\`{e}re energy for the fibrewise $S^{1}$-invariant plurisubharmonic functions in the "unit ball" under some assumptions. As an application we generalize the sharp Moser-Trudinger inequality on the unit ball.
\end{abstract}
\maketitle
\pagestyle{plain}

\section{Introduction}
In real analysis, given a real valued function $u$ defined on a domain $D$ in $\mathbb{R}^{n}$, its Schwarz symmetrization \cite{Kes06} is a function of the form $\widehat{u}=f(|x|)$, that is equidistributed with $u$.

One advantage of symmetrization is that we can reduce the study of many inequalities to the radial case, i.e. one variable case, for the reason that many quantities measuring the "size" of a function decrease under the symmetrization while any integrals of the form $\int_{D}F(u)dx$ stay. As an important application to complex geometry, it is discussed the behavior of a quantity called the Monge-Amp\`{e}re energy for a domain $D$ in $\mathbb{C}^{n}$ in \cite{Ber14} by Berman and Berndtsson. More specifically, if $u$ is a plurisubharmonic function in the unit ball vanishing on the boundary, it can be proved that the Monge-Amp\`{e}re energy decreases under the symmetrization provided that $u$ is $S^{1}$-invariant. Here the $S^{1}$-invariance of $u$ is assumed in order to maintain the plurisubharmonity under the symmetrization, which makes it meaningful to consider the energy after symmetrization.

In this article, we consider the domain in the total space $X$ of $L$ for a given negative line bundle on a compact complex manifold $Y$. It is a natural generalization of the case discussed in \cite{Ber14}. Indeed, $\mathbb{C}^{n}$ can be seen as the total space of the tautological line bundle $\mathcal{O}(-1)$ on $\mathbb{P}^{n-1}$, after blowing up at the origin. Then in order to define the symmetrization, we need a volume form (or measure, more analytically) $\Omega$ and a norm function $\Phi$ on $X$. In \cite{Ber14}, the norm function in $\mathbb{C}^{n}$ is chosen to be $|w|^{2}$, which can be seen as the push forward of the Fubini-Study metric on $\mathcal{O}_{\mathbb{P}^{n-1}}(-1)$:
$$|w_{1}|^{2}e^{\log(1+\sum_{i>1}|w_{i}/w_{1}|^{2})}$$
through the blow-up. Here we represent it in the coordinate patch $\{w_{1}\neq0\}$. Thus it is indicated in the model case that the functions with the form that $\Phi=|\xi|^{2}e^{\phi}$, where $\phi$ is a metric on $-L$ with positive curvature, are the suitable candidates if one wants the symmetrization to be related to the geometry of $Y$. Here we represent it in one local coordinate patch $(z,\xi)$.

We need some additional conditions about the volume form $\Omega$:
\begin{enumerate}
  \item In order to prove that the plurisubharmonicity is preserved after symmetrization, $\Omega$ needs to be homogeneous of some positive degree $2l$ with respect to $\xi$ and satisfy $\mathrm{Ric}(\Omega)\geqslant0$;

  \item In order to prove the symmetrization inequality of the Monge-Amp\`{e}re energy by variation argument, we need to pick a reference function with the form that $u_{0}=F(\Phi)$, which is plurisubharmonic on $X$, such that $u_{0}$ solves the following equation:
      $$\mathrm{MA}(u_{0})=G(\Phi)\Omega$$
      for some smooth function $G$.
\end{enumerate}
In order to figure out the meaning of these two requirements, we will do some simple calculation. Firstly we substitute $F(\Phi)$ to $\widetilde{F}(\log\Phi)$ (which means that $\widetilde{F}(t)=F(e^{t})$). Then out of the zero section,
$$\mathrm{MA}(u_{0})=(n+1)\widetilde{F}^{\prime\prime}(\widetilde{F}^{\prime})^{n}d\log\Phi\wedge d^{c}\log\Phi\wedge(dd^{c}\log\Phi)^{n},$$
which equals to
$$\frac{i}{\pi}(n+1)\widetilde{F}^{\prime\prime}(\widetilde{F}^{\prime})^{n}\omega_{\phi}^{n}\wedge\frac{d\xi\wedge d\overline{\xi}}{|\xi|^{2}}$$
$$=\frac{i}{\pi}(n+1)\Phi^{n+1}((F^{\prime})^{n}F^{\prime\prime}\Phi+(F^{\prime})^{n+1})\omega_{\phi}^{n}\wedge\frac{d\xi\wedge d\overline{\xi}}{|\xi|^{2}}.$$
Thus if we put
$$\widetilde{G}(t)=\frac{G(t)}{(n+1)((F^{\prime})^{n}(t)F^{\prime\prime}(t)t+(F^{\prime})^{n+1})(t)},$$
the second requirement implies that
$$\widetilde{G}(\Phi)\Omega=\frac{i}{\pi}\Phi^{n+1}\omega_{\phi}^{n}\wedge\frac{d\xi\wedge d\overline{\xi}}{|\xi|^{2}}.$$
Because $\Omega$ is assumed to be homogeneous of degree $2l$ with respect to $\xi$, $\widetilde{G}(t)$ must equal to $Ct^{n-l+1}$.

The only possibility is then that we consider the following volume form
$$\Omega=\frac{i}{\pi}\Phi^{l}\omega_{\phi}^{n}\wedge\frac{d\xi\wedge d\overline{\xi}}{|\xi|^{2}},$$
which is truly homogeneous of degree $2l$, and at this time the requirement about the Ricci curvature of $\Omega$ can be related to a curvature condition of $Y$:
$$\mathrm{Ric}(\omega_{\phi})\geqslant l\omega_{\phi}.$$
Unfortunately, they are not exactly equivalent because
$$\mathrm{Ric}(\Omega)=\mathrm{Ric}(\omega_{\phi})-l\omega_{\phi}-(l-1)[\xi],$$
Here $[\xi]$ refers to the current of integration of the zero section. Therefore $\mathrm{Ric}(\Omega)$ may not be semi-positive at the zero section if we only assume that $\mathrm{Ric}(\omega_{\phi})\geqslant l\omega_{\phi}$. Some problem does arise from this, and we will discuss them later. Anyway, the condition that
\[
\mathrm{Ric}(\omega_{\phi})\geqslant l\omega_{\phi}\tag{3}
\]
implies that $Y$ should be Fano and $lL-K_{Y}\geqslant0$. Thus concerning the Fano index (see \cite{Cas13} for the definition and some basic properties, for example), $l\leqslant n+1$ and the equality holds if and only if $(Y,L)=(\mathbb{P}^{n},\mathcal{O}(-1))$. Last but not least, it is worthwhile to point out that (3) is equivalent to the certain ("twisted") K\"{a}hler-Einstein equation \cite{Bm13,Don12,Sze11} on $Y$:
\[
\mathrm{Ric}(\omega_{\phi})=l\omega_{\phi}+\theta.
\]
Here $\theta$ is a fixed positive $(1,1)$-current, that may e.g., be the current of integration on a klt divisor. Therefore roughly speaking, they are equivalent to the existence of the certain ("twisted") K\"{a}hler-Einstein metric $\omega_{\phi}\in c_{1}(-L)$ on $Y$.

In summary, we should start with a negative line bundle $L$ on a Fano manifold $Y$ such that $lL-K_{Y}\geqslant0$, and assume that there is a metric $\phi$ on $-L$ such that $\omega_{\phi}:=dd^{c}\phi>0$ and
$$\mathrm{Ric}(\omega_{\phi})\geqslant l\omega_{\phi}.$$
Then we take the norm function on the total space $X$ of $L$ to be $\Phi:=|\xi|^{2}e^{\phi}$ and the volume form
$$\Omega_{l}=\frac{i}{\pi}\Phi^{l}\omega_{\phi}^{n}\wedge\frac{d\xi\wedge d\overline{\xi}}{|\xi|^{2}}.$$
If $u$ is a real valued function defined in a domain in $X$, we can define its Schwarz-type symmetrization to be the unique function of the form $\widehat{u}=f(\log\Phi)$, which is equidistributed with $u$ with respect to $\Omega$. Then the next step is to prove the plurisubharmonicity is preserved under symmetrization. Unfortunately, as is pointed before, the curvature condition
$$\mathrm{Ric}(\omega_{\phi})\geqslant l\omega_{\phi}$$
is not good enough. It cannot guarantee that $\mathrm{Ric}(\Omega_{l})\geqslant0$ at the zero section. However, the negative contribution of the term $|\xi|^{2}$ in $\Omega_{l}$ will always disappear when $l\leqslant1$, so we will consider these cases at first. We define the normalized the Monge-Amp\`{e}re energy by
$$\mathcal{E}(u)=\frac{1}{\mathrm{Vol}(-L)}\int_{\mathcal{B}}(-u)(dd^{c}u)^{n+1},$$
and our main result can be stated as following, which says that the main result in \cite{Ber14} extends perfectly to this more general case.

\begin{theorem}
Let $u$ be plurisubharmonic in the "unit ball" $\mathcal{B}:=\{\Phi<1\}$, and assume that $u$ extends continuously to the closed ball with zero boundary values. Assume also that $u$ is fibrewise $S^{1}$-invariant, and let $\widehat{u}$ be the Schwarz-type symmetrization of $u$. Then
$$\mathcal{E}(\widehat{u})\leqslant\mathcal{E}(u).$$
\end{theorem}

As an application, we generalize the sharp Moser-Trudinger inequality in the unit ball developed in \cite{Ber11}, to the "unit ball" defined above. The Moser-Trudinger inequality can be stated as follows:

\begin{theorem}
Let $u$ be a fibrewise $S^{1}$-invariant plurisubharmonic function in the "unit ball" in $X$ that vanishes on the boundary. Then if $\mathcal{E}(u)=1$,
$$\int_{\mathcal{B}}e^{l(-u)^{(n+2)/(n+1)}}\Omega_{l}\leqslant C/l$$
where $C$ is an absolute constant. In particular, this inequality implies the Moser-Trudinger inequality for $u$ with any bounded energy:
$$\log\int_{\mathcal{B}}e^{-u}\Omega_{l}\leqslant(\frac{n+1}{l})^{n+1}(\frac{1}{n+2})^{n+2}\mathcal{E}(u)+D$$
with $D=\log\frac{C}{l}$. Remember here $0<l\leqslant1$.
\end{theorem}

When $l>1$, it seems not to be totally hopeless. Indeed, when we look back to the model case in \cite{Ber14}, they didn't directly work on the total space $X$ of $L=\mathcal{O}(-1)$ but $X_{0}=\mathbb{C}^{n+1}$, which can be seen as the blow down of $X$ at the zero section. Because it is normal, the $L^{2}$-problem on $X_{0}$ is equivalent to be on $X$. The advantage of treating this ambient space instead is that the term $|\xi|^{2}$ won't negatively contribute neither, so the volume form $\Omega_{l}$ would be Ricci semi-positive for any $0<l\leqslant n+1$, after pulled back to the ambient space. It is indicated that for the general case, one could also consider the blow down at the zero section. The problem is that the ambient space may have singularity at the origin, except some special cases such as $(Y,L)=(\mathbb{P}(E),\mathcal{O}_{\mathbb{P}(E)}(-1))$ (which is pointed out by Mihai P\u{a}un) for some holomorphic vector bundle $E$. Probably our method could work in general, if one can analyze the behaviour of the Bergmann kernel (or holomorphic functions, more general) at the singular part. Currently, we would only work smoothly for $l>1$ provided the assumption that $X_{0}$ is smooth, and also get a family of Moser-Trudinger type inequalities as $l$ varies. Apparently in this occasion $X$ can be seen as $X_{0}$ blow up along some submanifold, and the upper bound of $l$ equals to the codimension of this submanifold, which we denote by $r$. Furthermore, we could take the coordinate patch of $X$ and $X_{0}$ to be $(x_{1},...,x_{r-1},z_{r},...,z_{n},\xi)$ and $(w_{1},...,w_{r},z_{r},...z_{n})$ respectively. In particular, if $(Y,L)=(\mathbb{P}(E),\mathcal{O}_{\mathbb{P}(E)}(-1))$, $r$ is just the rank of $E$ and $(w_{1},...,w_{r})$ the coordinate along fibre.

Similarly, we have the following inequalities:

\begin{theorem}
Let $X_{0}$ be the blow down of $X$ at the zero section. Assume that it is smooth. Let $u$ be an $S^{1}$-invariant with respect to $(w_{1},...,w_{r})$, plurisubharmonic function in the "unit ball" in the total space of $X_{0}$ that vanishes on the boundary. Then if $\mathcal{E}(u)=1$,
$$\int_{\mathcal{B}}e^{l(-u)^{(n+2)/(n+1)}}\Omega_{l}\leqslant C/l$$
where $C$ is an absolute constant. In particular, this inequality implies the Moser-Trudinger inequality for $u$ with any bounded energy
$$\log\int_{\mathcal{B}}e^{-u}\Omega_{l}\leqslant(\frac{n+1}{l})^{n+1}(\frac{1}{n+2})^{n+2}\mathcal{E}(u)+D$$
with $D=\log\frac{C}{l}$.
Remember here $0<l\leqslant r$.
\end{theorem}
These inequalities indicate that we can also do the Moser-Trudinger type estimate for some singular measure, i.e. $\Omega_{l}$, at the cost of the worse multiple constant $(\frac{n+1}{l})^{n+1}(\frac{1}{n+2})^{n+2}$ (comparing to the sharp estimate in \cite{Ber14}). In particular, when $(Y,L)=(\mathbb{P}^{n},\mathcal{O}(-1))$ and $\phi$ is the Fubini-Study metric, $\mathcal{B}$ will be the unit ball in the usual sense, and we can reformulate the Moser-Trudinger inequality as
$$\log\int_{B}e^{-u}\frac{d\lambda}{|w|^{2(n-l+1)}}\leqslant(\frac{n+1}{l})^{n+1}(\frac{1}{n+2})^{n+2}\mathcal{E}(u)+D.$$
Here $d\lambda$ is the Lebesgue measure. As we see, it returns back to the result in \cite{Ber14} when $l=n+1$.

\begin{acknowledgment}
The author would like to thank Professor Bo Berndtsson, who encourages the author to consider this problem and provides many useful suggestions.
\end{acknowledgment}

\section{Preliminaries}
Let $Y$ be an $n$-dimensional Fano manifold with a negative line bundle $L$. Assume that there is a metric $\phi$ on $-L$ such that $\omega_{\phi}:=i\partial\overline{\partial}\phi>0$ and
$$\mathrm{Ric}(\omega_{\phi})\geqslant l\omega_{\phi}$$
for some $0<l\leqslant1$. It means that there is a global non-negative function $\Phi$ on the total space $X$ of $L$, which locally can be written as
$$\Phi(z,\xi)=|\xi|^{2}e^{\phi(z)}.$$
Here $(z.\xi)$ is the local coordinate patch of $L$. We choose it as the norm function. Moreover, we define the volume form
$$\Omega_{l}:=\frac{i}{\pi}\Phi^{l}\omega_{\phi}^{n}\wedge\frac{d\xi\wedge d\overline{\xi}}{|\xi|^{2}}.$$
$\Omega_{l}$ can be seen as a measure of $X$. Then we come to our definition of Schwarz-type symmetrization on the total space of $L$.
\begin{definition}
If $u$ is a real valued function defined in a domain $D$ in $X$, its Schwarz-type symmetrization, is a fibrewise radial function,
$$\widehat{u}(z,\xi)=f(\log\Phi),$$
with $f$ increasing, that is equidistributed with $u$ with respect to $\Omega_{l}$. The latter requirement means that for any real $t$,
$$|\{u<t\}|_{\Omega_{l}}=|\{\widehat{u}<t\}|_{\Omega_{l}}=:\sigma(t).$$
\end{definition}

Before the further discussion, we recall some notions and properties about the Monge-Amp\`{e}re energy and geodesics needed later. One can refer to \cite{Ber11} for the full details.

Given a $u$ in PSH(D) we define its Monge-Amp\`{e}re energy by
$$\mathcal{E}(u)=\frac{1}{\mathrm{Vol}(-L)}\int_{D}(-u)(dd^{c}u)^{n+1}.$$
We consider the curve in the space of plurisubharmonic functions on $D$, which means the function
$$u_{t}(z,\xi)=u(t,z,\xi)$$
with a real parameter $t$, varying between $0$ and $1$. By definition, $u_{t}$ is a subgeodesic if $u(\textrm{Re}\tau,z,\xi)$ is plurisubharmonic as a function of $(\tau,z,\xi)$, and it is a geodesic if moreover this plurisubharmonic function solves the homogeneous complex Monge-Amp\`{e}re equation
$$(dd^{c}u_{\mathrm{Re}\tau})^{n+2}=0.$$
Here we can see $t$ as the real part of a complex variable $\tau\in\mathbb{C}$. We list some useful properties of $\mathcal{E}$ without proof.
\begin{proposition}

1.If $u_{t}$ is of class $C^{1}$, then $\mathcal{E}(u_{t})$ is differentiable with derivative
$$\frac{d}{dt}\mathcal{E}(u_{t})=\int_{D}-\frac{du_{t}}{dt}(dd^{c}_{z,\xi}u_{t})^{n+1};$$

2.If $u_{t}$ is moreover smooth,
$$dd^{c}_{t}\mathcal{E}(u_{t})=-\int_{D}(dd^{c}_{t,z,\xi}u_{t})^{n+2};$$

3.By the formula in 2 we conclude that $\mathcal{E}(u_{t})$ is an affine function of $t$ along any geodesic and concave along any subgeodesic. Furthermore we remark here that by a standard approximation technique the smoothness condition is not necessary.
\end{proposition}

Finally, we finish this section by two lemmas which will be used several times in the future. The first one is based on a result in \cite{Ber06}. We consider a pseudoconvex domain $\mathcal{D}$ in $\mathbb{C}\times X$ and its $(n+1)$-dimensional slices
$$D_{t}=\{(z,\xi)\in X;(t,z,\xi)\in\mathcal{D}\}$$
where $t$ ranges over (an interval in) $\mathbb{C}$. We say that a domain $D$ in $X$ is fibrewise $S^{1}$-invariant if $D$ is invariant under the map
$$(z,\xi)\mapsto(z,e^{i\theta}\xi)$$
for all $\theta$ in $\mathbb{R}$. A function $u$ is fibrewise $S^{1}$-invariant if $u(z,e^{i\theta}\xi)=u(z,\xi)$ for all real $\theta$.
\begin{lemma}
Assume that $\mathcal{D}$ is a pseudoconvex domain in $\mathbb{C}\times X$ such that all its slices $D_{t}$ are connected and fibrewise $S^{1}$-invariant. Assume also that the zero section belongs to $D_{t}$ when $t$ lies in a domain $U$ in $\mathbb{C}$. Then
$$\log|D_{t}|_{\Omega_{l}}$$
is a superharmonic function of $t$ in $U$.
\begin{proof}
First we claim that if $D$ is a fibrewise $S^{1}$-invariant domain that contains the zero section, then any fibrewise $S^{1}$-invariant holomorphic function $u$ on $D$ is constant. Indeed, $u$ can be locally expanded as
$$u(z,\xi)=\sum_{k=0}^{\infty} u_{k}(z)\xi^{k}.$$
Then $u(z,\xi)=u_{0}(z)$ by fibrewise $S^{1}$-invariance. On the other hand, restrict $u$ to the zero section $0$, we have
$$u(z,0)=u_{0}(z).$$
Since $u(z,0)$ is a global holomorphic function defined on $Y$, we conclude that $u_{0}(z)=u(z,0)$ is constant. We remark here that the same thing holds for the anti-holomorphic functions.

Now we consider the Bergman kernel $B_{t}((z,\xi),(w,\zeta))$ defined on $D_{t}$ with weight function given by $\Omega_{l}$. Fix a point $p\in Y$, then the Bergman kernel $B_{t}(p,(w,\zeta))$ is an anti-holomorphic function on $D_{t}$.  Since $D_{t}$ is fibrewise $S^{1}$-invariant by assumption, $B_{t}(p,(w,\zeta))$ is also fibrewise $S^{1}$-invariant. Hence
$$B_{t}(p,(w,\zeta))=B_{t}(p,p)$$
is constant, and since
$$\int_{D_{t}}1B_{t}(p,(w,\zeta))\Omega_{l}=1$$
we have
$$\log|D_{t}|_{\Omega_{l}}=-\log B_{t}(p,p).$$
Then the superharmonicity of $\log|D_{t}|_{\Omega_{l}}$ is a consequence of the main result in \cite{Ber06}, which says that $\log B_{t}(p,p)$ is subharmonic provided the weight function of $\Omega_{l}$ is plurisubharmonic.
\end{proof}
\end{lemma}
\begin{remark}
The situation here is not exactly the same as is in \cite{Ber06}. In fact, $D_{t}$ is assumed to be a domain in $X$ instead of in $\mathbb{C}^{n+1}$. But it won't make too much difference, since basically the proof is to use the reproducing property
$$B_{t}((z,\xi),(z,\xi))=\int_{D_{t}}B_{t}((w,\zeta),(z,\xi))\overline{B_{t}((w,\zeta),(z,\xi))}\Omega_{l}$$
to compute $\partial^{2}B_{t}((z,\xi),(z,\xi))/\partial t\partial\overline{t}$. The crucial step is that if we denote the weight function by $\psi$, then $u=\partial^{\psi}B_{t}$ is the minimal solution to the $\overline{\partial}$-equation
$$\overline{\partial}u=B_{t}\overline{\partial}\partial_{t}\psi.$$
Then by H\"{o}rmander's $L^{2}$-estimate for the holomorphic functions, we eventually get that
$$\frac{\partial^{2}B_{t}((z,\xi),(z,\xi))}{\partial t\partial\overline{t}}\geqslant\int_{D_{t}}|B_{t}|^{2}C\Omega_{l},$$
where
$$C=\psi_{t\overline{t}}-\sum(\psi_{(z,\xi)_{j}\overline{(z,\xi)_{k}}})^{-1}\psi_{t\overline{(z,\xi)_{j}}}\overline{\psi_{t\overline{(z,\xi)_{k}}}}.$$
Here $\psi_{(z,\xi)_{j}}$ means that we take derivative with respect to the $j$-th coordinate in $(z_{1},...,z_{n},\xi)$. Hence $C$ equals precisely the determinant of the full complex Hessian of $\psi$ divided by the determinant of the Hessian of $\psi_{t}$, which is positive. Certainly we do not definitely have the $L^{2}$-estimate for the holomorphic functions on the domain of the total space of $L$. However, we can always see a function $u$ as a $-K_{L}$-valued $(n+1,0)$-form. Then the thing also goes well because the assumption that $\mathrm{Ric}(\omega_{\phi})-l\omega_{\phi}\geqslant0$ provides the necessary positivity of the $L^{2}$-estimate for the holomorphic form, i.e. the volume form $\Omega_{l}$ defined before gives the metric on $-K_{L}$ with positive curvature.
\end{remark}

The next lemma could be well-known to experts. One could also refer to \cite{Ber14} for the proof.
\begin{lemma}
Let $u$ be a smooth function, then
$$D:=\{(w,z);u(z)-\mathrm{Re}w<0\}$$
is pseudoconvex iff $u$ is plurisubharmonic.
\end{lemma}

\section{Symmetrization of plurisubharmonic functions}
In order to discuss the symmetrization inequality on $X$, we need to verify that the symmetrization of a (fibrewise $S^{1}$-invariant) plurisubharmonic function is again plurisubharmonic. First we prove the following lemma.
\begin{lemma}
Let $u$ be a smooth plurisubharmonic function defined in an open set $U$ in $X$, and assume that $u$ vanishes on the boundary of $U$. Let
$$\sigma(t):=|\{(z,\xi);u(z,\xi)<t\}|_{\Omega_{l}}$$
for $t<0$. Then $\sigma$ is strictly increasing on the interval $(\min u,0)$.
\begin{proof}(sketch)
The proof is the same as the one of Lemma 2.2 in \cite{Ber14}. First, $\sigma(t)$ is certainly increasing. Thus if $\sigma(t)$ is constant in an interval $(t,t+\varepsilon)$, there would be some $s\in(t,t+\varepsilon)$ which is a regular value of $u$ by Sard's theorem. By the Hopf's lemma, the gradient of $u$ does not vanish on the boundary of $\{u<s\}$ unless $u$ is constant in $\{u<s\}$. In the latter case $s\leqslant\min u$. If this is not the case, i.e. $|\nabla u|>0$ on an open subset near the boundary of $\{u<s\}$, it is easy to see that $\sigma(t)$ must be strictly increasing at $s$.
\end{proof}
\end{lemma}

We can now prove the next result. We say that a domain $D$ in $X$ is balanced if for any $\lambda$ in $\mathbb{C}$ with $|\lambda|\leqslant1$ and $(z,\xi)$ in $D$, $(z,\lambda\xi)$ also lies in $D$.
\begin{theorem}
Let $D$ be a balanced domain in $X$ containing the zero section. Let $u$ be a fibrewise $S^{1}$-invariant plurisubharmonic function in $D$. Then $\widehat{u}$ is also plurisubharmonic.
\begin{proof}
We may of course assume that $u$ is smooth so that the previous lemma applies. By definition, $\widehat{u}$ can be written as
$$\widehat{u}=f(\log\Phi),$$
so what we need to prove is that $f$ is convex. Since for any real $t$,
$$\sigma(t):=|\{u<t\}|_{\Omega_{l}}=|\{\widehat{u}<t\}|_{\Omega_{l}}=|\{|\xi|^{2}e^{\phi}<e^{f^{-1}(t)}\}|_{\Omega_{l}}.$$
Hence
$$f^{-1}(t)=\frac{1}{l}\log\sigma(t)+C.$$
Here we use the fact that $\Omega_{l}$ is homogeneous of degree $2l$ with respect to $\xi$, hence $|\{|\xi|^{2}e^{\phi}<e^{f^{-1}(t)}\}|_{\Omega_{l}}$ is kind of like the volume of an $l$-dimensional ball $B_{l}$ with radius $|\xi|^{2}=e^{f^{-1}(t)-\phi}$. Since $\sigma$ is increasing, $f^{-1}$ is also increasing. Therefore $f$ is convex precisely when $f^{-1}$ is concave, i.e. when $\log\sigma$ is concave.

Consider the domain in $\mathbb{C}\times X$
$$\mathcal{D}=\{(\tau,z,\xi);(z,\xi)\in X, u-\mathrm{Re}\tau<0\}.$$
Then, if $t=\textrm{Re}\tau$, $\sigma=|D_{\tau}|_{\Omega_{l}}$. Note that $\mathcal{D}$ is pseudoconvex since $\phi-\textrm{Re}\tau$ is plurisubharmonic and we claim that $\mathcal{D}$ also satisfies all the other conditions of Lemma 2.1.

Let $(z,\xi)$ lie in $D_{\tau}$ for some $\tau$. The function $\gamma(\lambda):=u(z,\lambda\xi)$ is then subharmonic in the unit disk, and moreover it is radial by assumption that $u$ is fibrewise $S^{1}$-invariant, i.e.
$$\gamma(\lambda)=g(|\lambda|),$$
where $g$ is increasing. Therefore the whole disk $\{\lambda\xi\}$ is contained in $D_{\tau}$. In particular the zero section lies in any $D_{\tau}$. and $(z,0)$ can be connected with $(z,\xi)$ by a curve, so $D_{\tau}$ is connected. Thus Lemma 2.1 can be applied and we conclude that
$$\log\sigma(\textrm{Re}\tau)=\log|D_\tau|_{\Omega_{l}}$$
is a superharmonic function of $\tau$. Since this function only depends on $\textrm{Re}\tau$ it is actually concave, and the proof is complete.
\end{proof}
\end{theorem}

\section{The symmetrization inequality on "unit ball"}
In this section we focus on the main result of this article. In order to prove Theorem 1.1 we use a 2-variables version of the fact that the inverse of an increasing concave function is convex.
\begin{lemma}
Let $a(s,t)$ be a concave function of two real variables. Assume $a$ is strictly increasing with respect to $t$, and let $t=k(s,x)$ be the inverse of $a$ with respect to the second variable for $s$ fixed, so that $a(s,k(s,x))=x$. Then $k$ is convex as a function of both variables $s$ and $x$.
\end{lemma}

This lemma is sort of geometrically obvious, and one could refer to \cite{Ber14} for the proof.

Now we turn to the symmetrization inequality.
\begin{theorem}
Let $u$ be plurisubharmonic in the "unit ball" $\mathcal{B}:=\{\Phi<1\}$, and assume that $u$ extends continuously to the closed ball with zero boundary values. Assume also that $u$ is fibrewise $S^{1}$-invariant, and let $\widehat{u}$ be the Schwarz-type symmetrization of $u$. Then
$$\mathcal{E}(\widehat{u})\leqslant\mathcal{E}(u).$$
\begin{proof}
We do the variation along (sub)geodesics to confirm our result. Firstly, we need to prove that the symmetrization of an $S^{1}$-invariant (sub)geodesic is still a subgeodesic.
\begin{lemma}
Let $u_{t}$ be a subgeodesic of $S^{1}$-invariant plurisubharmonic functions. Then $\widehat{u}_{t}$ is also a subgeodesic.
\begin{proof}(of lemma)
Let $u_{s}$ be a subgeodesic which we may assume to be smooth. Let
$$A(s,t)=|\{u_{s}<t\}|_{\Omega_{l}}.$$
It follows again from Lemma 2.1 that $a:=\log A$ is a concave function of $s$ and $t$ together. As in the proof of Theorem 3.1 all we need to prove is that the inverse of $a$ with respect to $t$ is convex with respect to $s$ and $t$ jointly. But this is precisely the content of Lemma 4.1.
\end{proof}
\end{lemma}
The strategy to prove the symmetrization inequality is as follows. We put $u_{0}=F(\Phi)$, satisfying $F(1)=0$ and an equation
$$(dd^{c}u_{0})^{n+1}=G(u_{0})\Omega_{l},$$
where $G$ is some smooth function of a real variable. For example, we may take $u_{0}=\Phi-1$. We can also assume that $u_{1}=u$ is smooth by standard approximation procedure, then apply Chen's theorem in \cite{Che00} to construct a geodesic $u_{t}$ of class $C^{(1,1)}$ that connect $u_{0}$ and $u_{1}$. In fact, since we have assumed that $u_{0}$ and $u_{1}$ are smooth up to the boundary, we can by a max construction assume that they are both equal to $A\log((1+|\xi|^{2}e^{\phi})/2)$ for some large $A>0$, when $|\xi|^{2}e^{\phi}>1-\varepsilon$. Then $u_{0}$ and $u_{1}$ can be extended to plurisubharmonic functions in all of $X$, equal to $A\log((1+|\xi|^{2}e^{\phi})/2)$ outside of the unit ball. We can even consider them as the metrics on a line bundle $\mathcal{O}_{\mathbb{P}(L\oplus\mathcal{O}_{Y})}(A)$ over $\mathbb{P}(L\oplus\mathcal{O}_{Y})$. In fact, through the natural embedding
\begin{equation*}
\begin{split}
X&\hookrightarrow\mathbb{P}(L\oplus\mathcal{O}_{Y})\\
(z,\xi)&\mapsto(z,[\xi:1])=(z,[w_{1}:w_{2}])
\end{split}
\end{equation*}
$A\log((1+|\xi|^{2}e^{\phi})/2)$ can be seen as the restriction of a logarithmically homogeneous function $A(\log((|w_{2}|^{2}+|w_{1}|^{2}e^{\phi})/2)$ from the total space of $L\oplus\mathcal{O}_{Y}$ to $X$, which corresponds to a metric on $\mathcal{O}_{\mathbb{P}(L\oplus\mathcal{O}_{Y})}(A)$. Then as is stated in Chen's theorem, the space of K\"{a}hler potential on a compact K\"{a}hler manifold is geodesic convex of class $C^{(1,1)}$, which exactly gives the geodesics we need.

Now we consider the energy functionals along the two curves $u_{t}$ and $\widehat{u}_{t}$, $\mathcal{E}(u_{t})=:g(t)$ and $\mathcal{E}(\widehat{u}(t))=:h(t)$. Since $u_{0}$ is already 'symmetric', $g(0)=h(0)$, and we want to prove that $g(1)\geqslant h(1)$. We know that $g$ is affine and that $h$ is concave from Proposition 2.1 and Lemma 4.2, so this follows if we can prove that $g^{\prime}(0)=h^{\prime}(0)$. But
$$g^{\prime}(0)=\int-\frac{du_{0}}{dt}(dd^{c}u_{0})^{n+1},$$
by Proposition 2.1. We also claim that we can arrange things so that
$$h^{\prime}(0)=\int-\frac{d\widehat{u}_{t}}{dt}|_{t=0}(dd^{c}u_{0})^{n+1}.$$
By the choice of $u_{0}$,
$$g^{\prime}(0)=\int-\frac{du_{0}}{dt}G(u_{0})\Omega_{l}=\frac{d}{dt}|_{t=0}\int-H(u_{t})\Omega_{l},$$
if $H^{\prime}=G$. Similarly
$$h^{\prime}(0)=\frac{d}{dt}|_{t=0}\int-H(\widehat{u}_{t})\Omega_{l}.$$
But, since $u_{t}$ and $\widehat{u}_{t}$ are equidistributed
$$\int-H(u_{t})\Omega_{l}=\int-H(\widehat{u}_{t})\Omega_{l}$$
for all $t$. Hence $g^{\prime}(0)=h^{\prime}(0)$ and the proof is complete.

It remains to handle the claim about the derivative of $h$. $\widehat{u}_{0}=u_{0}$ is smooth and we can approximate $\widehat{u}_{1}$ from above by a smooth "symmetric" plurisubharmonic function. Now connect these two smooth functions by a geodesic, $v_{t}$, which can be taken to be $C^{(1,1)}$ by the argument above. Let
$$\mathcal{E}(v_{t})=:k(t).$$
Since $v_{t}\geqslant\widehat{u}_{t}$, $-\frac{d}{dt}v_{0}\leqslant-\frac{d}{dt}\widehat{u}_{0}$. We then apply the above argument to $k$ instead of $h$ and find that $k(1)\leqslant g(1)$. Taking limits as $v_{1}$ tends to $\widehat{u}_{1}$ we conclude the proof.
\end{proof}
\end{theorem}

\section{A Moser-Trudinger inequality for fibrewise $S^{1}$-invariant functions}
\subsection{The original version}
As an application, we can prove a Moser-Trudinger inequality for the fibrewise $S^{1}$-invariant functions on $\mathcal{B}=\{\Phi<1\}$.
\begin{theorem}
Let $u$ be a fibrewise $S^{1}$-invariant plurisubharmonic function in the "unit ball" in $X$ that vanishes on the boundary. Then if $\mathcal{E}(u)=1$
$$\int_{\mathcal{B}}e^{l(-u)^{(n+2)/(n+1)}}\Omega_{l}\leqslant C/l$$
where $C$ is an absolute constant. In particular, this inequality implies the Moser-Trudinger inequality for $u$ with any bounded energy:
$$\log\int_{\mathcal{B}}e^{-u}\Omega_{l}\leqslant(\frac{n+1}{l})^{n+1}(\frac{1}{n+2})^{n+2}\mathcal{E}(u)+D$$
with $D=\log\frac{C}{l}$. Remember here $0<l\leqslant1$.
\begin{proof}
We need the following result of Moser \cite{Mos71}.
\begin{lemma}(Moser)
If $w$ is an increasing function on $(-\infty,0)$ that vanishes when $t$ goes to zero and satisfies
$$\int^{0}_{-\infty}(-w^{\prime})^{n+2}dt\leqslant1$$
then
$$\int^{0}_{-\infty}e^{(-w)^{(n+2)/(n+1)}}e^{t}dt\leqslant C,$$
where $C$ is an absolute constant.
\end{lemma}
First, we scale the Moser's inequality stated in Lemma 5.1. Applying Lemma 5.1 to $w_{k}(s):=k^{(n+1)/(n+2)}w(s/k)$ we obtain that
$$\int^{0}_{-\infty}e^{k(-w)^{(n+2)/(n+1)}}e^{kt}dt\leqslant C/k,$$
under the same hypothesis. Now in order to prove our conclusion, we can assume that $u=f(\log\Phi)$ by symmetrization inequality. Moreover, we can assume that $f(t)$ is constant for $t$ sufficiently large negative, and the general by approximation. The advantage of this assumption is that the behaviour on the boundary would be better at this case, then we can freely use the Fubini's theorem. We can even assume that $f$ is smooth. Indeed, if we have proved the estimate above for the smooth ones, then for a general $u$ with $\mathcal{E}(u)=1$, we can always approximate it from above by a sequence of smooth functions $\{u_{t}\}$ with $\mathcal{E}(u_{t})$ tending to $\mathcal{E}(u)$, and apply the estimate to $u_{t}/\mathcal{E}(u_{t})^{\frac{1}{n+2}}$. The estimate for the general one then follows by taking the limit. Now let
$$F(t)=e^{l(-t)^{(n+2)/(n+1)}},$$
then
\begin{equation*}
\begin{split}
&\int_{\mathcal{B}}e^{l(-u)^{(n+2)/(n+1)}}\Omega_{l}=\int_{\mathcal{B}}F(u)\Omega_{l}\\
                                                                            &=\frac{1}{\pi}\int_{Y}e^{l\phi}\omega_{\phi}^{n}\int_{|\xi|^{2}<e^{-\phi}} F\circ f(\log\Phi)|\xi|^{2(l-1)}id\xi\wedge d\overline{\xi}\\
                                                                            &=\frac{1}{\pi}\int_{Y}e^{l\phi}\omega_{\phi}^{n}\int_{t^{2}<e^{-\phi}}F\circ f(\log t^{2}+\phi)t^{2l}\frac{dt\wedge d\theta}{t}\\
                                                                            &=\int_{Y}e^{l\phi}\omega_{\phi}^{n}\int^{-\phi}F\circ f(s+\phi)e^{ls}ds\\
                                                                            &=\int_{Y}e^{l\phi}\omega_{\phi}^{n}\int^{0}F\circ f(t)e^{l(t-\phi)}dt\\
                                                                            &=\mathrm{Vol}(-L)\int^{0}F\circ f(t)e^{lt}dt
\end{split}
\end{equation*}
We substituted several times during the computation. Also notice that here $C^{\prime}:=\mathrm{Vol}(-L)$ is a constant only depends on $L$.

In order to apply the Moser's inequality, we need to estimate the term $\int^{0}_{-\infty}(f^{\prime})^{n+2}dt$. Indeed, we could do the following calculation: if we write $\log\Phi$ as
$$\log\Phi=\log|\xi|^{2}+\phi,$$
then outside the zero section, we have
\begin{equation*}
\begin{split}
(dd^{c}u)^{n+1}&=(n+1)f^{\prime\prime}(f^{\prime})^{n}d\log\Phi\wedge d^{c}\log\Phi\wedge(dd^{c}\log\Phi)^{n}\\
               &=(n+1)f^{\prime\prime}(f^{\prime})^{n}d\log\Phi\wedge d^{c}\log\Phi\wedge(dd^{c}\phi)^{n}
\end{split}
\end{equation*}
Thus by Fubini's theorem,
\begin{equation*}
\begin{split}
\mathcal{E}(u)&=\frac{n+1}{\mathrm{Vol}(-L)}\int_{Y}(dd^{c}\phi)^{n}\int^{0}-ff^{\prime\prime}(f^{\prime})^{n}(t)dt\int_{\log\Phi=t}d^{c}\log\Phi\\
              &=\frac{n+1}{\mathrm{Vol}(-L)}\int_{Y}(dd^{c}\phi)^{n}\int^{0}-ff^{\prime\prime}(f^{\prime})^{n}(t)dt\int_{\log\Phi<t}dd^{c}\log\Phi
\end{split}
\end{equation*}

Since $\log\Phi$ is fibrewise logarithmically homogeneous,
$$dd^{c}\log\Phi(z,\lambda\xi)=dd^{c}\log\Phi(z,\xi)$$
while $dd^{c}_{\xi}\log\Phi(z,\lambda\xi)=\lambda^{2}(dd^{c}_{\xi}\log\Phi(z,\xi))$ by direct calculation, which means that it is harmonic with respect to $\xi$ outside of the zero section. But $\log\Phi-\log|\xi|^{2}$ is bounded near the zero section so this point mass must be the same as
$$dd^{c}\log|\xi|^{2}.$$
Hence
\begin{equation*}
\begin{split}
\mathcal{E}(u)&=(n+1)\int^{0}-ff^{\prime\prime}(f^{\prime})^{n}(t)dt\\
              &=\int^{0}(f^{\prime})^{n+2}(t)dt
\end{split}
\end{equation*}
which shows the conclusion. The general case, when $f$ is not constant near $-\infty$ follows from approximation. All in all, $\int^{0}_{-\infty}(f^{\prime})^{n+2}dt$ equals to the Monge-Amp\`{e}re energy of $u$, which is $1$.

Therefore we can apply the scaled version of Lemma 5.1 with $k=l$ and
$$w=f.$$
The first inequality follows.

For the latter part, we start with the elementary inequality for positive numbers $x$ and $y$
$$xy\leqslant\frac{1}{n+2}x^{n+2}+\frac{n+1}{n+2}y^{(n+2)/(n+1)}.$$
This implies
$$y\leqslant\frac{1}{n+2}x^{n+1}+\frac{n+1}{n+2}y^{(n+2)/(n+1)}/x.$$
Choose $x$ so that
$$x^{n+1}=(\frac{n+1}{l(n+2)})^{n+1}\mathcal{E}(u)$$
and take $y=(-u)$. Then
\begin{equation*}
\begin{split}
&-u\leqslant(\frac{n+1}{l})^{n+1}\frac{1}{(n+2)^{n+2}}\mathcal{E}(u)+l(-\frac{u}{\mathcal{E}(u)^{1/(n+2)}})^{(n+2)/(n+1)}
\end{split}
\end{equation*}
Therefore the former inequality implies the Moser-Trudinger inequality for fibrewise $S^{1}$-invariant functions
$$\log\int e^{-u}\Omega_{l}\leqslant(\frac{n+1}{l})^{n+1}(\frac{1}{n+2})^{n+2}\mathcal{E}(u)+D.$$
Obviously we apply the former inequality to $\frac{u}{\mathcal{E}(u)^{1/(n+2)}}$ here.
\end{proof}
\end{theorem}

\subsection{A special case}
It is believed that the sharp multiplicity constant of the Moser-Trudinger inequality
$$\log\int_{\mathcal{B}}e^{-u}\leqslant C\mathcal{E}(u)+D$$
for a pseudoconvex domain $\mathcal{B}$ in $\mathbb{C}^{n+1}$ should be $C=\frac{1}{(n+2)^{n+2}}$. However, inspired by the Moser-Trudinger inequality proved before, it is possible to consider some singular measures instead of the Lebesgue measure.

First, we start with a measure $e^{-\chi}d\mu$ on $\mathbb{C}^{n+1}$, where $d\mu$ is the Lebesgue measure and $\chi$ is assumed to be log-homogenous of degree $2(n-l+1)$. We can pull it back through the blow-up
$$f:\mathcal{O}(-1)\rightarrow\mathbb{C}^{n+1}.$$
If we focus on one coordinate ball $U_{1}=\{w_{1}\neq0\}$, and denote the coordinate of $\mathcal{O}(-1)$ by $(\xi,z_{i})$, it can be written as
\begin{equation*}
\begin{split}
&e^{-\chi(\xi,\xi z_{i})}d\xi\wedge d\overline{\xi}\wedge d(\xi z_{i})\wedge d(\overline{\xi z_{i}})\\
=&|\xi|^{2n}e^{-\chi(\xi,\xi z_{i})}d\xi\wedge d\overline{\xi}\wedge dz_{i}\wedge d\overline{z_{i}}\\
=&|\xi|^{2n}|\xi|^{-2(n-l+1)}e^{-\chi(z)}d\xi\wedge d\overline{\xi}\wedge dz\wedge d\overline{z}\\
=&|\xi|^{2l}e^{-\chi}\frac{d\xi\wedge d\overline{\xi}}{|\xi|^{2}}\wedge dz\wedge d\overline{z}\\
=&(|\xi|^{2}e^{\phi})^{l}e^{-(\chi+l\phi)}\frac{d\xi\wedge d\overline{\xi}}{|\xi|^{2}}\wedge dz\wedge d\overline{z}
\end{split}
\end{equation*}
Here $\phi$ is an auxiliary metric fixed before. Thus this measure, after pull back, can be seen as a measure on the total space of $\mathcal{O}(-1)$, which is homogeneous of degree $2l$ with respect to $\xi$, and the assumption the weight function being plurisubharmonic (or using the former language that the Ricci curvature of this volume form is semi-positive) is satisfied as long as $l\leqslant n+1$.

On the other hand, given a measure $e^{-\eta}\frac{d\xi\wedge d\overline{\xi}}{|\xi|^{2}}\wedge dz\wedge d\overline{z}$ on the total space of $\mathcal{O}(-1)$, we can push it forward after multiplying with $(|\xi|^{2}e^{\phi})^{l}$. Note that at this time it is homogeneous of degree $2l$ with respect to the variable $\xi$.
\begin{equation*}
\begin{split}
&(|w_{1}|^{2}e^{\phi(w_{i}/w_{1})})^{l}e^{-\eta(w_{i}/w_{1})}\frac{dw_{1}\wedge d\overline{w_{1}}}{|w_{1}|^{2}}\wedge d(w_{i}/w_{1})\wedge d\overline{w_{i}/w_{1}}\\
=&|w_{1}|^{2l}|w_{1}|^{-2(n+1)}e^{-(\eta-\phi)(w_{i}/w_{1})}dw\wedge d\overline{w}\\
=&e^{-(\eta-\phi)(w_{i}/w_{1})-2(n-l+1)\log|w_{1}|}dw\wedge d\overline{w}
\end{split}
\end{equation*}
It is a measure on $\mathbb{C}^{n+1}$ whose weight function is log-homogenous of degree $2(n-l+1)$. Thus there is a correspondence between the measure whose weight function is log-homogenous of degree $2(n-l+1)$ on $\mathbb{C}^{n+1}$ and the homogeneous of degree $2l$ measure on the total space of $\mathcal{O}(-1)$. Observe that the negative contribution of the term $|\xi|^{2}$ (or $|w_{1}|^{2}$ after push forward) is absorbed by the basis $dw\wedge d\overline{w}$, we can get rid of the restriction of $l$ before, i.e. $0<l\leqslant1$, to consider any $l$ such that $0<l\leqslant n+1$ instead.

Inspired by this observation, one can even consider a more general case, namely the total space of $L$ which is smooth after blow down at the zero section. We denote the blow down by $X_{0}$. We believe that the smoothness is not necessary, but it would require something like $L^{2}$-estimate on the singular variety, which is not so clearly currently.

We start from the same place as before, i.e. $(X,\Phi,\Omega_{l})$ with the certain conditions. We consider the blow down $X_{0}$ of $X$. If the local coordinate patch of $X$ is taken to be $(x_{1},...,x_{r-1},z_{r},...,z_{n},\xi)$, where $\{x_{1},...,x_{r-1},\xi\}$ are the $r$ coordinates involved in the blow down, $\Phi$ locally can be written as
$$\Phi(x,z,\xi)=|\xi|^{2}e^{\phi(x,z)}.$$
and the volume form
$$\Omega_{l}:=\Phi^{l}\omega_{\phi}^{n}\wedge\frac{d\xi\wedge d\overline{\xi}}{|\xi|^{2}}.$$
We push them forward to the ambient space $X_{0}$, and still denote them by $\Phi$ and $\Omega_{l}$. Moreover if we denote the coordinate of $X_{0}$ by $(w_{1},...,w_{r},z_{r},...,z_{n})$, and focus on one coordinate patch $\{w_{1}\neq0\}$, $\Omega_{l}$ is a measure on the total space of $X_{0}$ with a singular weight function
$$2(r-l)\log|w_{1}|+\psi(w_{i}/w_{1},z).$$
Here $\psi$ is a plurisubharmonic function given by $e^{l\phi}\omega_{\phi}^{n}$ after push forward. Next we define the Schwarz-type symmetrization on the total space of $X_{0}$ with respect to $\Phi$ and
$$\Omega_{l}=(\frac{1}{|w_{1}|^{2}})^{r-l}e^{-\psi(w_{i}/w_{1},z)}dw\wedge d\overline{w}\wedge dz\wedge d\overline{z}$$
as usual. Remember here $(\frac{1}{|w_{1}|^{2}})^{r-l}e^{-\psi(w_{i}/w_{1},z)}$ is the representative of a global function (with respect to $(w_{1},...,w_{r})$) in one coordinate patch, $\{w_{1}\neq0\}$. Then it turns out that this symmetrization keeps the plurisubharmonicity of the fibrewise $S^{1}$-invariant function defined on the balanced domain of the total space of $X_{0}$ by the same argument as is in Theorem 3.1. Certainly the fibre here refers to $(w_{1},...,w_{r})$. Moreover, the symmetrization inequality is still valid under this circumstance. Notice that basically the reason that the symmetrization inequality works is that there is a reference function $u_{0}=F(\log\Phi)$ solves the following Monge-Amp\`{e}re equation:
$$(dd^{c}u_{0})^{n+1}=G(u_{0})\Omega_{l}$$
with some smooth function of a real variable $G$. However this thing is guaranteed through the construction of $\Phi$ and $\Omega_{l}$ as is explained in Introduction. Since the symmetrization inequality is confirmed, we can use the usual method to prove the Moser-Trudinger inequality.
\begin{theorem}
Let $u$ be a fibrewise $S^{1}$-invariant plurisubharmonic function in the "unit ball" defined by $\mathcal{B}:=\{\Phi<1\}$ in the total space of $X_{0}$ that vanishes on the boundary. Then if $\mathcal{E}(u)=1$
$$\int_{\mathcal{B}}e^{l(-u)^{(n+2)/(n+1)}}\Omega_{l}\leqslant C/l$$
where $C$ is an absolute constant. In particular, this inequality implies the Moser-Trudinger inequality for $u$ with any bounded energy:
$$\log\int_{\mathcal{B}}e^{-u}\Omega_{l}\leqslant(\frac{n+1}{l})^{n+1}(\frac{1}{n+2})^{n+2}\mathcal{E}(u)+D$$
with $D=\log\frac{C}{l}$. Remember here $0<l\leqslant r$.
\begin{proof}
The proof is the same as of Theorem 5.1.
\end{proof}
\end{theorem}

Note that the "unit ball" may be fruitful as $l$ varies, which provides a lot of pseudoconvex domains in the total space of $X_{0}$. In particular, if $r=n+1$, and take $Y$ to be the projective space with $\phi$ the Fubini-Study metric, then $\mathcal{B}$ would be the unit ball in the usual sense. Therefore we get a family of Moser-Trudinger inequalities on the unit ball mentioned before.
\begin{corollary}
Let $u$ be an $S^{1}$-invariant plurisubharmonic function in the unit ball that vanishes on the boundary. Then if $\mathcal{E}(u)=1$
$$\int_{B}e^{l(-u)^{(n+2)/(n+1)}}\frac{d\lambda}{|w|^{2(n-l+1)}}\leqslant C/l$$
where $C$ is an absolute constant. In particular, this inequality implies the Moser-Trudinger inequality for $u$ with any bounded energy:
$$\log\int_{B}e^{-u}\frac{d\lambda}{|w|^{2(n-l+1)}}\leqslant(\frac{n+1}{l})^{n+1}(\frac{1}{n+2})^{n+2}\mathcal{E}(u)+D$$
with Lebesgue measure $d\lambda$ and $D$ an absolute constant. Remember here $0<l\leqslant n+1$.
\end{corollary}

Observe that when $l=n+1$, it returns back to the result in \cite{Ber11}, but the general one cannot be easily deduced from the one in \cite{Ber11} (at least it's not an obvious consequence). Indeed, one may try to use H\"{o}lder's inequality to approach it. Taking the conjugate exponents $\frac{n+1}{l},\frac{n+1}{n-l+1}$, it can be calculated as follows:
\begin{equation*}
\begin{split}
         &\int_{B}e^{l(-u)^{(n+2)/(n+1)}}\frac{d\lambda}{|w|^{2(n-l+1)}}\\
\leqslant&(\int_{B}e^{(n+1)(-u)^{(n+2)/(n+1)}}d\lambda)^{\frac{l}{n+1}}(\int_{B}\frac{d\lambda}{|w|^{2(n+1)}})^{\frac{n-l+1}{n+1}}
\end{split}
\end{equation*}
However, the term $\int_{B}\frac{d\lambda}{|w|^{2(n+1)}}$ is unbounded, thus we will get nothing. On the other hand, the multiple constant $(\frac{n+1}{l})^{n+1}(\frac{1}{n+2})^{n+2}$ is sort of sharp, because $\int_{B}\frac{d\lambda}{|w|^{2(n+1-\varepsilon)}}$ is definitely bounded.

\address{

\small Current address: Mathematical Sciences-Chalmers University of Technology, SE-412 96, Gothenburg, Sweden

\small E-mail address: jincao@chalmers.se
}


\begin{thebibliography}{99}
\bibitem{Don12}
S. Donaldson, K\"{a}hler metrics with cone singularities along a divisor, Essays in mathematics and its applications, 49-79.
\end{thebibliography}
\end{document}